\documentclass[12pt]{article}
\usepackage[a4paper,left=2cm,right=1.2cm,top=3cm,bottom=4cm]
{geometry}

\usepackage{amsmath,amstext, amsthm, empheq, extpfeil,  
amsbsy,amssymb,marvosym,fancyhdr,graphicx,amscd,amsfonts,latexsym,delarray,stackrel,lineno,color,cite,appendix,xcolor,url,hyperref}
\usepackage{wasysym}
\usepackage{mdframed}

\usepackage{srcltx}
\usepackage{mathrsfs}
\usepackage{float} 
\usepackage{subfig} 
\usepackage[all]{xy}
\usepackage{t1enc}
\usepackage{mathrsfs}
\usepackage{pifont}

\usepackage{mathpple}
\usepackage[T1]{fontenc}

\newcommand*\patchAmsMathEnvironmentForLineno[1]{%
  \expandafter\let\csname old#1\expandafter\endcsname\csname #1\endcsname
  \expandafter\let\csname oldend#1\expandafter\endcsname\csname end#1\endcsname
  \renewenvironment{#1}%
     {\linenomath\csname old#1\endcsname}%
     {\csname oldend#1\endcsname\endlinenomath}}%
\newcommand*\patchBothAmsMathEnvironmentsForLineno[1]{%
  \patchAmsMathEnvironmentForLineno{#1}%
  \patchAmsMathEnvironmentForLineno{#1*}}%
\AtBeginDocument{%
\patchBothAmsMathEnvironmentsForLineno{equation}%
\patchBothAmsMathEnvironmentsForLineno{align}%
\patchBothAmsMathEnvironmentsForLineno{flalign}%
\patchBothAmsMathEnvironmentsForLineno{alignat}%
\patchBothAmsMathEnvironmentsForLineno{gather}%
\patchBothAmsMathEnvironmentsForLineno{multline}%
}

\definecolor{Green}{rgb}{0,1,0}
\definecolor{Blue}{RGB}{0,0,191}
\definecolor{mathmodecolor}{RGB}{0,102,0}
\definecolor{keywordcolor}{RGB}{0,51,151}
\definecolor{sourcebackgroundcolor}{RGB}{255,247,223}
\definecolor{unixagred}{RGB}{255,0,0}
\definecolor{lightgray}{RGB}{191,191,191}
\definecolor{green}{RGB}{1,191,191}

\newtheorem{thm}{Theorem}[section]

\newtheorem{cor}[thm]{Corollary}
\newtheorem{lem}[thm]{Lemma}

\newtheorem{defn}[thm]{Definition}

\newtheorem{quest}[thm]{Question}

\def\F{{\mathbb F}}

\def\qqq{\,,\quad \forall}

\def\bm2{{\rm \B mod^2}}
\def\b2{{\rm \B mod^{\mathfrak s}}}

\def\tr{{\rm tr}}

\def\C{{\mathbb C}}
\def\F{{\mathbb F}}

\def\N{{\mathbb N}}

\def\Q{{\mathbb Q}}
\def\R{{\mathbb R}}

\def\Z{{\mathbb Z}}

\def\B{{\mathbb B}}

\def\tr{{\rm tr}}

\def\cR{{\mathcal R}}

\def\part{\partial}

\newcommand{\ie}{{\it i.e.\/}\ }

\definecolor{trust}{rgb}{0,1,1}

\def\noncommutative geometry{{noncommutative geometry }}

\def\qqq{\,,\quad \forall}

\title
{The  Carlitz group of the rationals}

\begin{document}

\maketitle

\centerline{Alain Connes} 

\vspace{1cm}
\centerline{\em{To the memory of David Goss, a great  number theorist friend.}}

\vspace{1cm}

\begin{abstract}
This paper contains two parts. The first is the solution  of a challenge question, proposed by Etienne Ghys, on the determination of all maps from rational numbers to themselves such that the difference quotient (f(x)-f(y))/(x-y) is always a square. The second is the  computer determination, done with the help of Stephane Gaubert, of a function of primes which plays a key role  in the first part as a generalization of a result of Carlitz.
\end{abstract}
\vspace{0.1in}

\section{Introduction}

Let $K$ be a field and $G(K)$ be the group, which we call the Carlitz group of $K$,  of all bijections $f:K\to K$ which fulfill the following condition 
\begin{equation}\label{sq condition0}
\frac{f(x)-f(y)}{x-y}\in K^2 \qqq x\neq y \in K	
\end{equation}
These bijections form a group under composition. This group is the group of all bijections when $K$ is algebraically closed or when it is a perfect field of characteristic two. When $K=\R$ it is the group of orientation preserving homeomorphisms of the line. For finite fields of odd characteristic it was determined by Carlitz \cite{Carlitz} as the semi-direct product of the affine group by the Frobenius automorphisms, a group already considered by Galois in his work on primitive solvable equations. 

We consider the following  very intriguing question formulated by Etienne Ghys\footnote{as a challenge during a meeting of the French Academy of Sciences} 
\begin{quest} (E. Ghys)
Determine $G(\Q)$.	
\end{quest}

Our main result is to determine all  maps $f:\Q\to \Q$ such that the following holds 
\begin{equation}\label{sq condition}
\frac{f(x)-f(y)}{x-y}\in \Q^2 \qqq x\neq y	
\end{equation}

The answer is given by the following Theorem:

\begin{thm}\label{main} A map $f:\Q\to \Q$,  fulfills 
	\eqref{sq condition} if and only if it is is an affine map, $f(x)=a^2 x+b$.
\end{thm}
In order to prove Theorem \ref{main} we first refine the result of Carlitz for finite fields of odd characteristic by showing (Theorem \ref{arbp} below) that a self-map of $\F_q$ which fulfills \eqref{sq condition0} is either constant or bijective. As a corollary of this preliminary result we get the {\em simplicity} of the Paley graphs. The notion of simplicity for a graph is straightforward (see Definition \ref{gammacong} below). We then focus in \S \ref{ratsequences} on sequences of rational numbers and show that if $f:\N\to \Q$ is such that 	\eqref{sq condition} holds for $x\neq y\in \N$ and that $f(0)=0$, $f(1)=1$, then $f(j)=j$ for all $j\in \N$. With this at hand one easily gets the Theorem \ref{main}. 

The second part of the paper (\S \ref{sectlp}) is based on extensive computations done in collaboration with Stephane Gaubert and we are indebted to him for his great help. 
It deals with  a quantitative form of Theorem \ref{arbp}.   Our result on infinite sequences of rational numbers does not exclude the existence of arbitrarily long finite sequences of rational numbers which fulfill the same conditions. This leads one to study for each (odd) prime $p$  the function $L(p)$ which associates to the prime $p$ the smallest number $L$ such that $f(j)=j$ is the only solution of  $$f(0)=0, \ f(1)=1, \ \frac{f(x)-f(y)}{x-y}\in \F_p^2 \qqq x\neq y	\in \{0,\ldots, L\}$$
In fact we consider the function of two variables $W(p,L)$ (where $p$ is an odd prime and $ L<p$ an integer), which gives the number of solutions of the above equation. We give a simple Gaussian estimate of $W(p,L)$ in \S \ref{gaussect}. 
We then  show in \S \ref{sectlowbound} that the function $L(p)$ is larger than the function $n(p)$ giving the first quadratic non-residue. In particular known results on the latter show that the growth of $L(p)$ cannot be $O(\log(p))$ in spite of the slow growth of the function $L(p)$ for primes up to $p=443$. 
The computation in \S \ref{exptests} of  the function $W(p,x)$ for sufficiently many primes, in order to make educated guesses on its general behavior, was done by heavy use of parallel computations\footnote{The author thanks the M\'esocentre Phymath federating the CMAP, CMLS, CPHT and PMC laboratories of \'Ecole polytechnique,  for providing access to the cluster ``Hopper'' were parallel computation were performed. Thanks also to François Bachelier and to Pieter van Bijnen for their help at an earlier stage of the computations.}. The results show that the Gaussian estimate is good in many cases but one meets several primes for which the computation of $W(p,L)$ and of $L(p)$ requires much longer than what the Gaussian estimate would suggest. We measure the non-gaussian behavior of such primes by the function $\sigma(p)=\sum \log W(p,k)$ and  compare it with the behavior in $(\log p)^3$ given by the Gaussian estimate.

\section{The Carlitz group of $\Q$}
We give in this section the proof of Theorem \ref{main}. We first refine, in \S \ref{strongcarlitz}, the result of Carlitz to remove the hypothesis of bijectivity. We then apply this result to sequences of rational numbers in \S \ref{ratsequences} and complete the proof of Theorem \ref{main} in \S \ref{sectmain}.
\subsection{Strengthening of the result of Carlitz}\label{strongcarlitz}

 Let  $p$ be an odd prime, $q$ a power of $p$ and $\chi:\F_q\to \{-1,0,1\}\subset \C$ denote the quadratic residue character\footnote{$\chi:\F_q\to \C$, $\chi(a)=0\iff a=0$, $\chi(ab)=\chi(a)\chi(b)$, $a\in \F_q^2\iff \chi(a)\neq -1$. When $q=p$ is prime it is the Legendre symbol.}.
	  \begin{lem}\label{gammap} Let $H\subset F_q$ be a subset of cardinality $>1$ and with non-empty complement.\newline
	  $(i)$~Assume that $\chi(x-y)=\chi(x'-y)$ for all $x,x'\in H$ and $y\notin H$. Then the  cardinality of $H$ fulfills $\#H\leq (q-1)/2$.\newline
	  $(ii)$~Assume that $\#H\leq (q-1)/2$ then for any pair of distinct elements $u,v\in H$ there exists at least two elements $y\notin H$ such that $\chi(u-y)\neq \chi(v-y)$.\newline
	  $(iii)$~There exists $u,v\in H$, $y\notin H$ such that $\chi(u-y)\neq \chi(v-y)$.
	   \end{lem}
	   \proof $(i)$~By a translation one can assume that $0\notin H$. Then $\chi(x)=\chi(x')$ for all $x,x'\in H$, and thus $H$ is contained in one of the halves of the multiplicative group given by squares or non-squares. \newline
	  $(ii)$~Assume that $\#H\leq (q-1)/2$ then its complement $H^c$ contains at least $(q-1)/2+1$ elements. Moreover one has the classical formula\footnote{with $\chi$ taking values in $\C$} (\cite{Jones}, Lemma 2.1)
	  $$
	  \sum_{z\in\F_q}\chi(u-z)\chi(v-z)=-1
	  $$
	  The number $j$ of terms equal to $-1$ in the sum over $z\in H$ is  at most $\#H-2$  since the contributions of $z=u$ and $z=v$ vanish. Let $k$ be the number of terms equal to $-1$ in the sum over $z\in H^c$. Then one has $j+k$ terms equal to $-1$, two equal to zero and $q-(j+k)-2$ terms equal to $1$ thus
	  $$
	  -(j+k)+(q-(j+k)-2)=-1\Rightarrow j+k=(q-1)/2
	  $$
	  But one has $j\leq \#H-2$ and $\#H\leq (q-1)/2$ thus one gets $k\geq 2$.\newline
	  $(iii)$~Assume that $\chi(x-y)=\chi(x'-y)$ for all $x,x'\in H$ and $y\notin H$. Then by $(i)$, the cardinality of $H$ fulfills $\#H\leq (q-1)/2$. Thus $(ii)$ applies and for any pair of distinct elements $u,v\in H$ there exists at least two elements $y\notin H$ such that $\chi(u-y)\neq \chi(v-y)$. Thus one gets a contradiction. \endproof 
	   \begin{thm}\label{arbp} Let $q$ be a power of an odd prime and $f:\F_q\to \F_q$ such that the following holds 
\begin{equation}\label{sq condition1bis}
 \frac{f(x)-f(y)}{x-y}\in \F_q^2 \qqq x\neq y	
\end{equation}
then $f$ is  an affine map times a power of the Frobenius automorphism. It is either constant or bijective.
\end{thm}
\proof  It is enough to show that if $f$ is not constant it is injective since then the result is Carlitz's Theorem \cite{Carlitz}. Assume that $H$ is a non-trivial fiber of $f$ so that 
$$
H=\{u\in \F_q\mid f(u)=a\}, \ \ \#H\geq 2, \ \#H^c>0
$$
Let us show that $\chi(x-y)=\chi(x'-y)$ for all $x,x'\in H$ and $y\notin H$. One has
$$
f(x)=f(x'), \ f(x)-f(y)\neq 0, \ f(x')-f(y)\neq 0
$$	
and thus \eqref{sq condition1bis} in the form  
$$
\left(\frac{f(x)-f(y)}{x-y}\right)\left(\frac{f(x')-f(y)}{x'-y}\right)\in \F_q^2
$$   
implies
$$
(x-y)(x'-y)\in \F_q^2
$$
which in turns means that $\chi(x-y)=\chi(x'-y)$. We can thus apply Lemma \ref{gammap} and get a contradiction. \endproof

To state the geometric corollary of Theorem \ref{arbp} for the Paley graphs, we introduce the following notion of simplicity for graphs:
 \begin{defn}\label{gammacong}
	 $(i)$~An equivalence relation $\cR$ on the vertices of a graph $\Gamma$ is a $\Gamma$-congruence 	if and only if for two distinct $\cR$-classes $C,C'$ the fact that $(x,x')$ is or is not an edge is independent of the choices of   $x\in C$ and $x'\in C'$.\newline
	 $(ii)$~A graph $\Gamma$ is {\em simple} if and only if the only $\Gamma$-congruence are the two trivial ones\footnote{the diagonal and the coarse one}.
	 \end{defn}
	 When $q$ is a power of an odd prime, and is congruent to $1$ modulo $4$ one defines the Paley graph $\Gamma(q)$ as the graph with set of vertices $V=\F_q$ and where two vertices $x,y$ are adjacent if and only if $\chi(x-y)=1$. 

\begin{cor} The Paley graphs are simple.	
\end{cor}
\proof Let $\cR$ be a $\Gamma(q)$-congruence. Choose a section $C\mapsto a(C)\in C$ \ie an element in each equivalence class. Let $f:\F_q\to\F_q$ be the projection $f(x):=a(C(x))$. Let us show that \eqref{sq condition1bis} holds. For $x,x'\in \F_q$, either $f(x)=f(x')$ and \eqref{sq condition1bis} holds or the $\cR$-classes $C,C'$ of $x,x'$ are distinct. In that case the fact that $(x,x')$ is or is not an edge is independent of the choices of   $x\in C$ and $x'\in C'$ and we can thus choose the elements $a(C)$ and $a(C')$. This shows that $\chi(x-x')=\chi(f(x)-f(x'))$ and hence that \eqref{sq condition1bis} holds.
Applying Theorem \ref{arbp} to $f$ one gets the simplicity of $\Gamma(q)$.\endproof

\subsection{Sequences of rationals}\label{ratsequences}

 The problem for $\Q$ gives the following question : study all sequences $f(j)\in \Q$, $j\in \N$ such that 
 \begin{equation}\label{sequcondition}
 \frac{f(i)-f(j)}{i-j}\in \Q^2 \qqq i\neq j	
\end{equation}
 Let $\rho_p:\Z_{(p)}\to \F_p$ be the morphism from the ring $\Z_{(p)}$ of fractions with denominator prime to $p$ to the quotient by the ideal generated by $p$.
 
\begin{lem}\label{lemloc} Let $p$ be a prime. Let $f(j)\in \Q$, $j\in \{0,\ldots ,p-1\}$ such that \eqref{sequcondition} holds for all pairs $i\neq j\in \{0,\ldots ,p-1\}$.  Assume  that the denominators of the $f(j)$ are not divisible by $p$ for $j<p$. 	Then the map $j\mapsto \rho_p(f(j))\in \F_p$, $j<p$,  is a self-map of $\F_p$ which fulfills \eqref{sq condition1bis}.
\end{lem}
\proof Since $f(j)\in \Z_{(p)}$ for $j<p$ one gets for $i,j\in \{0,\ldots ,p-1\}$ that $\frac{f(i)-f(j)}{i-j}\in \Z_{(p)}$. An element $z \in \Z_{(p)}$ which is a square in $\Q$ is a square in $\Z_{(p)}$ because squaring doubles the $p$-adic valuations, and one gets that the $p$-adic valuation of the rational square root of $z$ is $\geq 0$. It follows that $\frac{f(i)-f(j)}{i-j}\in \Z_{(p)}^2$ and one gets, applying $\rho_p$ that the map $j\mapsto \rho_p(f(j))\in \F_p$ fulfills \eqref{sq condition1bis}.\endproof

\begin{lem}\label{propconj}
Let $f:\N\to \Z$ be such that 	\eqref{sequcondition} holds and that $f(0)=0$, $f(1)=1$. Then $f(j)=j$ for all $j\in \N$.
\end{lem}
\proof For each prime $p$ one can consider the map $f_p:\F_p\to \F_p$ obtained by reducing $f(j)$ modulo $p$ for $0\leq j<p$. For two distinct elements $0\leq i< j<p$ one has by hypothesis that $(j-i)(f(j)-f(i))$ is the square of a rational number, and hence the square of an integer. Thus the same holds for $f_p$ and by Theorem \ref{arbp}, one has that $f_p$ is the identity. This shows that  $f(j)-j$ is divisible by any prime $>j$ and hence is equal to $0$. \endproof

\begin{lem}\label{propconj1}
Let $f:\N\to \Q$ be such that 	\eqref{sequcondition} holds and that $f(0)=0$, $f(1)=1$. Then $f(j)=j$ for all $j\in \N$.
\end{lem}
\proof Let us  look at the denominators which may appear in a sequence $f(j)$ of rational numbers ($j\geq 0$) such  that $f(0)=0$, $f(1)=1$ and that \eqref{sequcondition} holds. For a prime $p$ we look at the first occurrence of a negative power of $p$ in $f(j)$:
\begin{equation}\label{jp def}
j_p(f):=\inf\{j\mid {\rm Val}_p(f(j))<0\}	
\end{equation}
One has $j_p(f)\geq 2$ by construction. Moreover since ${\rm Val}_p(f(j_p-1))\geq 0$ one has 
$$
{\rm Val}_p(f(j_p))={\rm Val}_p(f(j_p)-f(j_p-1))\in 2\Z
$$
since $f(j_p)-f(j_p-1)$ is a square. In fact one can consider the finite differences
$$
\alpha(k)=(f(j_p+k)-f(j_p-1))/(1+k), \ \ k\in \{0,\ldots, p-2\}
$$
By hypothesis $\alpha(k)$ are squares in $\Q$ and since $(1+k)$ is prime to $p$,  the $p$-adic valuation  ${\rm Val}_p(\alpha(k))$, if it is negative, is the same as ${\rm Val}_p(f(j_p+k)-f(j_p-1))$. It is even since $\alpha(k)$ are squares and thus if $k\in \{0,\ldots, p-2\}$ is such that  
${\rm Val}_p(f(j_p+k))<0$ one gets 
$$
{\rm Val}_p(f(j_p+k))={\rm Val}_p(f(j_p+k)-f(j_p-1)))\in 2\Z  
$$
We now consider the interval of length $p$ given by $I_p:=\{j_p-2,\ldots , j_p+p-3\}$ and we assume $p\geq 3$ so that $I_p$ contains $j_p$. We then define 
\begin{equation}\label{ep def}
e_p(f):=\inf\{{\rm Val}_p(f(j))\mid j\in I_p\}	
\end{equation}
Since all the numbers ${\rm Val}_p(f(j))$, $j\in I_p$, which are negative are even we get that $e_p(f)$ is even and $<0$.
We can thus multiply the $f(j)$ for $j\in I_p$ by $p^{-e_p(f)}$ without altering the fact that the finite differences are squares of rationals. We then consider the map
\begin{equation}\label{ep use}
r_p(f)(k):=\rho_p(f(j_p-2+k)p^{-e_p(f)})\in \F_p	
\end{equation}
The map $r_p(f)$ fulfills \eqref{sq condition1bis} by  Lemma \ref{lemloc}. It takes the same value $0$ at $k\in \{0,1\}$. This implies, by Theorem \ref{arbp} that it is constant equal to $0$ but this gives a contradiction since there exists a non-zero value of 
$r_p(f)$ due to the definition of $e_p(f)$. Thus we cannot have a non-trivial denominator involving odd primes. 

Let us now consider the case $p=2$. The definition \eqref{jp def} gives an integer $j_2$ and we assume $j_2<\infty$. One has ${\rm Val}_2(f(j_2))<0$. By construction one has $j_2\geq 2$ and
$$
{\rm Val}_2(f(j_2-2))\geq 0, \ \ {\rm Val}_2(f(j_2-1))\geq 0
$$
Thus, since $f(j_2)-f(j_2-1)$ is a square, 
$$
{\rm Val}_2(f(j_2))={\rm Val}_2(f(j_2)-f(j_2-1))\in 2\Z
$$
But $(f(j_2)-f(j_2-2))/2$ is a square and this gives a contradiction since
$$
{\rm Val}_2((f(j_2)-f(j_2-2))/2)={\rm Val}_2(f(j_2)-f(j_2-2))-1={\rm Val}_2(f(j_2))-1
\in 1+ 2\Z
$$
We have shown that no denominator can appear in the sequence $f(j)$ and thus it is integer valued. But then we can apply Lemma \ref{propconj} to get the conclusion. \endproof 
\subsection{Proof of Theorem \ref{main}}\label{sectmain}
 By Lemma \ref{propconj1}
we know that the only sequences $a(n)$, $n\in \N$, of rational numbers, $a(0)=0$, such that 	\eqref{sequcondition} holds are in fact constant times $n$  (and the constant is equal to a square). Indeed either $a(n)$ is constant equal to $0$ or there exists a smallest $j_0>0$ for which $a(j_0)\neq 0$, but then $a(j_0)$ is a square since $a(j_0-1)=0$ and \eqref{sequcondition} holds. Thus  the sequence $b(u):=a(j_0+u-1)/a(j_0)$ fulfills the hypothesis of Lemma \ref{propconj1} and is hence equal to $u$ for all $u>0$. This gives $a(j_0+u-1)=u a(j_0)$ for all $u\geq 0$. Let us show that $j_0=1$. One has $b(u)=0$ for $u\in \{-j_0+1,\ldots,0\}$ and $b(1)=1$. Thus if $j_0\geq 2$, \eqref{sequcondition} gives that $(b(1)-b(-1))/2$ is a square which is a contradiction. Thus $j_0=1$, $a(u)=ua(1)$ for all $u\geq 0$. Let then  $f:\Q\to \Q$ which fulfills 
	\eqref{sq condition}.  One can assume that $f(0)=0$ by subtracting $f(0)$. Let $x\in \Q$, $x\neq 0$. Let $a(n):=f(n x) /x$. It fulfills the condition \eqref{sequcondition}. Thus we get $a(n)=na(1)$ for all $n\in \N$. This shows that  
	 $f(nx)/x=n f(x)/x$ and thus $f(nx)=n f(x)$ for all $n\in \N$. Thus we have for integers  $a,b>0$
	 $$
	 f(a/b)b=f(a)=af(1)\Rightarrow f(a/b)=a/b f(1) \Rightarrow f(x)=xf(1)\qqq x\in \Q_+.
	 $$
	 This result applies also to the function $g(x):=f(1)-f(1-x)$ and gives 
	 $$
	 g(x)=xg(1), \forall x\in \Q_+ \Rightarrow f(1)-f(1-x) =x f(1), \forall x\in \Q_+
	 \Rightarrow f(1-x)=(1-x)f(1), \forall x\in \Q_+.
	 $$
	 This shows that $f(x)=xf(1)$ for all $x\in \Q$. \endproof

	  \section{The minimal length $L(p)$}\label{sectlp}
	  
	  Theorem \ref{main} does not exclude the existence of finite sequences $f(j)\in \Q$, $j\leq L$ such that \eqref{sequcondition} holds for all pairs $i\neq j\in \{0,\ldots, L\}$.	  For instance one has the following examples of length $4$ and $5$, \ie $L=3$ and $L=4$,
	  $$
	\left\{0,1,\frac{15842}{1681},23763\right\}, \ \   \left\{0,1,\frac{2738}{2209},\frac{3267}{2209},\frac{5476}{2209}\right\}
	  $$
	 The above proof of Theorem \ref{main} suggests to first control, given a prime $p$, the possible sequences 
	 $f(j)\in \F_p$, $0\leq j\leq L<p$ such that 
	 	  \begin{equation}\label{sq condition1}
f(0)=0, \ f(1)=1, \ \frac{f(x)-f(y)}{x-y}\in \F_p^2 \qqq x\neq y	\in \{0,\ldots, L\}
\end{equation}
By Theorem \ref{arbp} there exists a smallest $L=L(p)<p$ such that the only solution of \eqref{sq condition1} is $f(j)=j$. 	We show  the graph of the  function of $n$ giving the minimal  length $L(p(n))$ with $p=p(n)$ the n-th prime. It is shown until $n=79$ \ie $p(n)=401$.
	  
	  \begin{figure}[H]
\begin{center}
\includegraphics[scale=0.7]{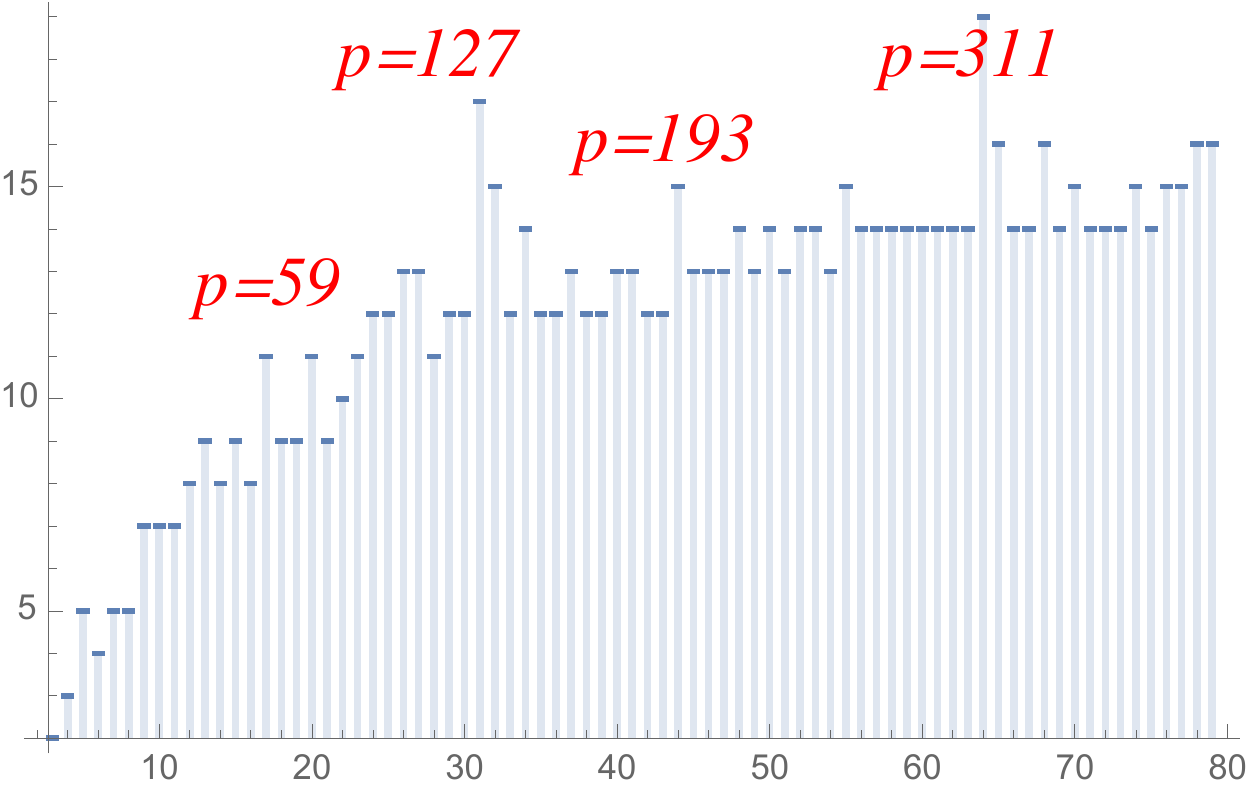}
\end{center}
\caption{Graph of the function $L(p(n))$ for $5\leq n\leq 79$. \label{wait1} }
\end{figure}

\subsection{Gaussian estimate of the function $W(p,L)$}\label{gaussect}
In order to give an estimate of the order of magnitude of the function $L(p)$ we do a simple counting. In fact we consider the function of two variables $W(p,L)$ where $p$ is an odd prime and $2\leq L<p$ an integer, which gives the number of solutions of \eqref{sq condition1}. Each time we write that some quantity is a square in $\F_p$ the probability that this is true is the fraction 
$$
P=\frac{1+(p-1)/2}{p}=\frac{p+1}{2p}
$$
When we deal with a list of the form 
$$
(0,1,f(2),\ldots , f(x))
$$
the number of variables is $x-1$ and the freedom is thus of $p^{x-1}$. The number of requirements that some expression is a square is 
$$
S(x)=\sum_{r=2}^x r=\frac{1}{2} \left(x^2+x-2\right)
$$
Thus if we assume the independence at the probabilistic level of the conditions we estimate the number of remaining possibilities as
$$
R=P^{S(x)} p^{x-1}. 
$$
This gives as an approximation  for $W(p,L)$ a Gaussian function $G(p,L)=e^{Q(p,L)}$ with exponent the quadratic form, in terms of $L=x$, 
$$
Q(p,x)=-\frac{1}{2} x^2 \log \left(\frac{2p}{p+1}\right)+\frac{1}{2} x \log \left(\frac{p(p+1)}{2}\right)-\log \left(\frac{p+1}{2}\right)
$$
This quadratic form vanishes for $x=1$ and the other root is thus 
$$
\ell(p):=2\,\frac{\log(p+1)-\log(2)}{\log(2)-\log(1+1/p)}\sim 
2\,\frac{\log(p)}{\log(2)}
$$
The next lemma shows that the probabilistic estimate is good at the beginning. Indeed for sequences 	$(0,1,f(2))$ there are two conditions and they should select about $1/4$ of the possible values of $f(2)\in \F_p$. More precisely one gets:
\begin{lem} Let $p$ be an odd prime. 
The number of sequences 	$(0,1,f(2))$ which fulfill the two conditions $f(2)/2\in \F_p^2$ and
$f(2)-1\in \F_p^2$ is equal to $E(p/4)+1$ where $E(x)$ denotes the integral part of $x$.
\end{lem}
\proof We consider the plane curve $C$ defined by the equation $1+x^2=2y^2$. By intersecting the lines through the point $(1,1)$ we get the rational parametrization: $t\mapsto P(t),$
\begin{equation}\label{curveC}
P(t):=(x,y)=\left(1-\frac{2 (2 t-1)}{2 t^2-1},1-\frac{2 t (2 t-1)}{2 t^2-1}\right),  \ t= (y-1)/(x-1)	
\end{equation}
For $x=1$ we have the two points with $y=\pm 1$. For $x\neq 1$ the ratio $t= (y-1)/(x-1)$ is finite and it determines uniquely $t$ and the point $P(t)$. The values of $t$ obtained from $x\neq 1$ are such that $t\neq 1/2$ and $t^2\neq 1/2$. Moreover for $x=-1$ one gets $t=0$ and $y=1$, or $t=1$ and $y=-1$. But all the 4 points $(\pm 1,\pm 1)$ correspond to the same solution $f(2)=2$.  Thus for the other solutions they correspond to values 
$$
A_p:=\{t\in \F_p \mid  t\neq 1/2, t^2\neq 1/2, t\neq 0, t\neq 1\}
$$
The cardinality of $A_p$ is $p-3$ if $2$ is not a square in $\F_p$ and is $p-5$ otherwise.

The two transformations $\alpha$, $(x,y)\mapsto (-x,y)$ and $\beta$, $(x,y)\mapsto (x,-y)$ can now be expressed as the following projective  involutions in the variable $t$, they correspond to the matrices
$$
\alpha=\left(
\begin{array}{cc}
 2 & -1 \\
 2 & -2 \\
\end{array}
\right), \ \ \beta=\left(
\begin{array}{cc}
 1 & -1 \\
 2 & -1 \\
\end{array}
\right), \ \alpha(t)=\frac{2 t-1}{2 (t-1)}, \  \beta(t)=\frac{t-1}{2 t-1}
$$
The fixed points of $\alpha$ are $\left\{t\to \frac{1}{2} \left(2-\sqrt{2}\right)\right\},\left\{t\to \frac{1}{2} \left(\sqrt{2}+2\right)\right\}$ which is empty unless $2$ is a square. The fixed points of $\beta$ are $\left\{t\to \frac{1}{2}-\frac{i}{2}\right\},\left\{t\to \frac{1}{2}+\frac{i}{2}\right\}$ which is empty unless $-1$ is a square. We need to consider 4 cases determined by the residue of $p$ modulo $8$.
\begin{itemize}
\item $p\equiv 1$ modulo $8$. Then both $-1$ and $2$ are squares. Both 
$\alpha$ and $\beta$ have two fixed points and thus the number of orbits of the group $H=\Z/4\Z$ acting on $A_p$ is 
$$
\# A_p/4+1=\frac{p-5}{4}+1=\frac{p-1}{4}=E(\frac p4)
$$
Thus adding the contribution of $(\pm 1,\pm 1)$ one gets the expected result.
\item $p\equiv 3$ modulo $8$. Then both $-1$ and $2$ are not squares. Both 
$\alpha$ and $\beta$ have no fixed points and thus the number of orbits of the group $H=\Z/4\Z$ acting on $A_p$ is 
$$
\# A_p/4=\frac{p-3}{4}=E(\frac p4)
$$
Thus adding the contribution of $(\pm 1,\pm 1)$ one gets the expected result.
\item $p\equiv 5$ modulo $8$. Then  $-1$ is a square and $2$ is not a square. Thus  
$\alpha$ has no fixed point but $\beta$ has two fixed points and thus the number of orbits of the group $H=\Z/4\Z$ acting on $A_p$ is 
$$
(\# A_p-2)/4 +1=\frac{p-5}{4}+1=E(\frac p4)
$$
Thus adding the contribution of $(\pm 1,\pm 1)$ one gets the expected result.
\item $p\equiv 7$ modulo $8$. Then  $-1$ is not a square and $2$ is  a square. Thus  
$\alpha$ has two fixed points but $\beta$ has no fixed points and thus the number of orbits of the group $H=\Z/4\Z$ acting on $A_p$ is 
$$
(\# A_p-2)/4 +1=\frac{p-7}{4}+1=E(\frac p4)
$$
Thus adding the contribution of $(\pm 1,\pm 1)$ one gets the expected result.	
\end{itemize}
This shows that in all cases the answer is $E(\frac p4)+1$. \endproof 

\begin{figure}[H]
\begin{center}
\includegraphics[scale=0.6]{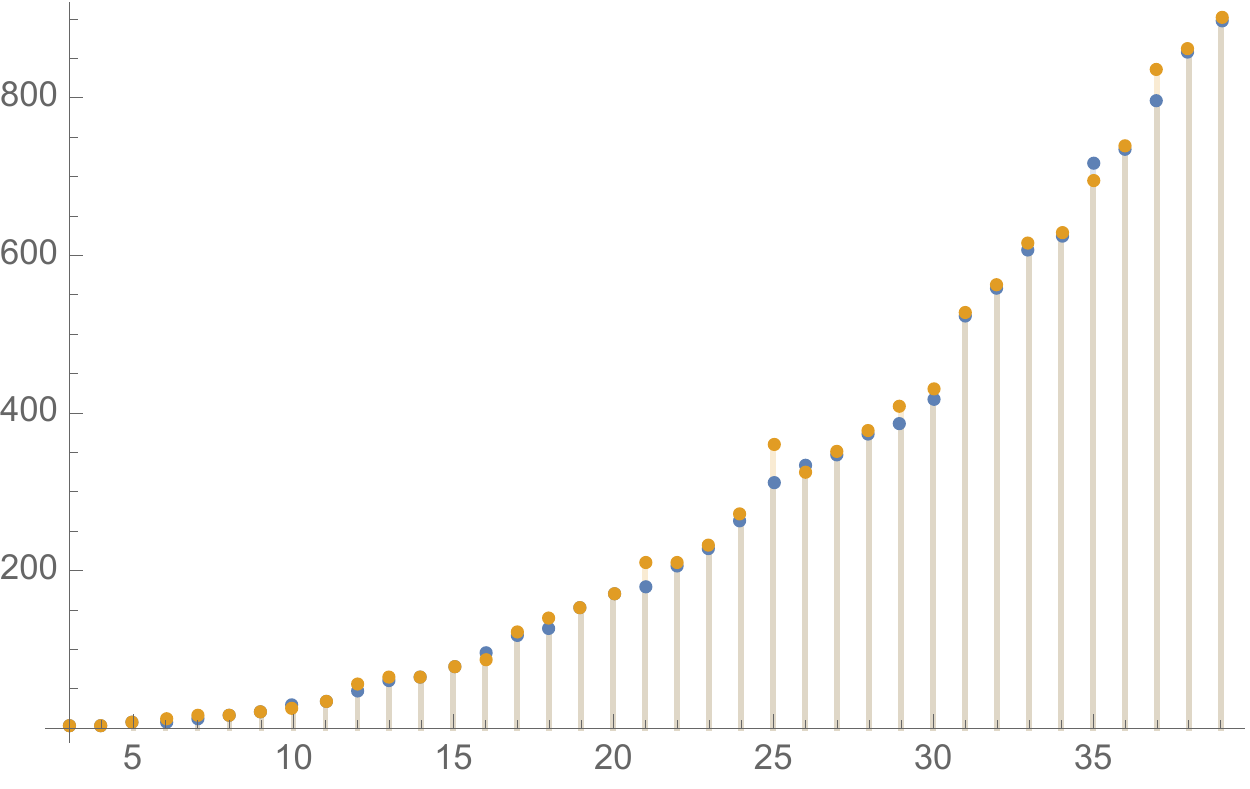}
\end{center}
\caption{Number of solutions compared to the Gaussian. \label{gaussiantolist2} }
\end{figure}

When we move to the next step \ie sequences of the form $
(0,1,f(2), f(3))
$, the probabilistic estimate is still rather good as shown in Figure \ref{gaussiantolist2}
In fact the discrepancy is governed by the family of elliptic curves $E(s)$, $v^2=Q_s(u)$ where 
$$
Q_s(u):=4 \left(4 s^4+16 s^3-28 s^2+8 s+1\right)u^4
-48 \left(2 s^2-1\right)^2 u^3$$ $$
+12 \left(36 s^4-16 s^3-12 s^2-8 s+9\right) u^2
-72 \left(2 s^2-1\right)^2 u
+9 \left(4 s^4+16 s^3-28 s^2+8 s+1\right)
$$
whose discriminant is 
$$
D(s)=2^{20}3^6\left(2 s^2-1\right)^4 \left(2 s^2-4 s+1\right)^4 \left(2 s^2-2 s+1\right)^4
$$
and does not vanish for rational values of $s$.
\begin{lem} Let $P(s)=(x,y)$ be a point of the curve $C$ of \eqref{curveC} defined by the equation $1+x^2=2y^2$.  
The  sequences 	$(0,1,f(2),f(3))$ which fulfill $f(2)=1+x^2=2y^2$ and  the  condition 
$
(f(i)-f(j))/(i-j)
$ is a square for $i<j$, $i,j\in \{0,\ldots ,3\}$
correspond to the points of the elliptic curve $E(s)$.
\end{lem}
\proof We consider the two conditions $f(3)-f(0)=3 X^2$ and $f(3)-f(1)=2 Y^2$. This gives the plane curve $C'$ defined by the equation $1+2 Y^2=3X^2$. By intersecting the lines through the point $(1,1)$ we get the rational parametrization: $u\mapsto R(u),$
\begin{equation}\label{curveC1}
R(u):=(X,Y), \ \ X=1-\frac{2 (2 u-3)}{2 u^2-3}, \ \ Y=\frac{2 u^2-6 u+3}{3-2 u^2}	
\end{equation}
One then writes the missing equation \ie that $
f(3)-f(2)
$ is a square. This gives the equation $v^2=Q_s(u)$ where the polynomial $Q_s(u)$ is determined by the equality 
$$
f(3)-f(2)=\frac{Q_s(u)}{\left(2 s^2-1\right)^2 \left(2 u^2-3\right)^2}
$$
using \eqref{curveC} to input $f(2)$ and \eqref{curveC1} to input $f(3)$ as 
$$
f(2)=\left(1-\frac{2 (2 s-1)}{2 s^2-1}\right)^2+1, \ \ f(3)=3 \left(1-\frac{2 (2 u-3)}{2 u^2-3}\right)^2
$$
Note that in the counting of points of the elliptic curve $E(s)$, the point at $\infty$ counts $2$ because it is a double point. This means that the number of points which are not at $\infty$ is given by $p-\tr F -1$ where $\tr F$ is the trace of the Frobenius.

\subsection{Lower bound for $L(p)$}\label{sectlowbound}

Let us define the function $n(p)$ for a prime $p$ as the largest integer $0\leq n<p$ such that 
all elements of $\{0,\ldots, n-1\}$ are quadratic residues. In other words $n(p)$ is the first quadratic non-residue. 
\begin{lem}\label{lowerbound} Let $p$ be an odd prime. One has $n(p)\leq L(p)$.
\end{lem}
\proof Consider the sequence given by 
$$
f(0)=0, \ f(u)=1\qqq u\in \{1,\ldots,n(p)-1\}
$$
Let us show that $f$ fulfills \eqref{sq condition1} for $u,v\in \{1,\ldots,n(p)-1\}$. One can assume that $v=0$ and $u\neq 0$ since if both $u,v$ are $\neq 0$ \eqref{sq condition1} is fulfilled since $f(u)=f(v)$. Then $f(u)-f(v)=1$ and \eqref{sq condition1} means that $u$ is a quadratic residue, which is true by definition of $n(p)$. 
Thus as long as $n(p)>2$ one has a sequence $f(i)$ which fulfills \eqref{sq condition1} for $i\leq n(p)-1$ and is not $f(i)=i$, thus $L(p)>n(p)-1$ and one gets $n(p)\leq L(p)$.\endproof


\begin{table}
\setlength{\tabcolsep}{3pt}
\begin{center}
\resizebox{16cm}{!}
{\begin{tabular}{cc|ccccccccccccccccccccccccc}
\fontfamily{ptm}\selectfont\scriptsize
&$p$&2&3&4&5&6&7&8&9&10&11&12&13&14&15&16&17&18&19&20&21\\
\hline
3&5&2&1&1&1&1&1&1&1&1&1&1&1&1&1&1&1&1&1&1&1\\
4&7&2&3&1&1&1&1&1&1&1&1&1&1&1&1&1&1&1&1&1&1\\
5&11&3&6&10&3&1&1&1&1&1&1&1&1&1&1&1&1&1&1&1&1\\
6&13&4&11&17&1&1&1&1&1&1&1&1&1&1&1&1&1&1&1&1&1\\
7&17&5&15&29&15&1&1&1&1&1&1&1&1&1&1&1&1&1&1&1&1\\
8&19&5&15&33&12&1&1&1&1&1&1&1&1&1&1&1&1&1&1&1&1\\
9&23&6&21&52&90&42&7&1&1&1&1&1&1&1&1&1&1&1&1&1&1\\
10&29&8&25&45&75&117&10&1&1&1&1&1&1&1&1&1&1&1&1&1&1\\
11&31&8&36&56&56&28&3&1&1&1&1&1&1&1&1&1&1&1&1&1&1\\
12&37&10&55&155&85&57&18&3&1&1&1&1&1&1&1&1&1&1&1&1&1\\
13&41&11&63&207&529&1529&616&92&4&1&1&1&1&1&1&1&1&1&1&1&1\\
14&43&11&66&248&690&788&166&2&1&1&1&1&1&1&1&1&1&1&1&1&1\\
15&47&12&78&312&813&847&266&12&3&1&1&1&1&1&1&1&1&1&1&1&1\\
16&53&14&85&265&575&1017&258&3&1&1&1&1&1&1&1&1&1&1&1&1&1\\
17&59&15&120&574&1050&810&409&107&32&18&9&1&1&1&1&1&1&1&1&1&1\\
18&61&16&139&621&1655&4109&2115&290&3&1&1&1&1&1&1&1&1&1&1&1&1\\
19&67&17&153&807&2908&4326&1685&99&9&1&1&1&1&1&1&1&1&1&1&1&1\\
20&71&18&171&944&3515&10090&18831&13816&3820&208&5&1&1&1&1&1&1&1&1&1&1\\
21&73&19&209&1651&6117&14375&6791&3795&1919&1&1&1&1&1&1&1&1&1&1&1&1\\
22&79&20&210&954&2368&2932&1676&360&26&6&1&1&1&1&1&1&1&1&1&1&1\\
23&83&21&231&1451&6030&10700&12441&13179&18342&23716&379&1&1&1&1&1&1&1&1&1&1\\
24&89&23&273&1707&7183&28187&25364&11066&774&12&3&2&1&1&1&1&1&1&1&1&1\\
25&97&25&359&3685&15195&37653&22445&10709&7201&8&3&2&1&1&1&1&1&1&1&1&1\\
26&101&26&325&2385&10255&30265&70983&144977&269175&343141&7300&161&4&1&1&1&1&1&1&1&1\\
27&103&26&351&2121&7077&12131&20426&31580&46109&64592&5120&56&3&1&1&1&1&1&1&1&1\\
28&107&27&378&2946&14571&31389&23486&3784&684&59&3&1&1&1&1&1&1&1&1&1&1\\
29&109&28&407&2969&11767&36697&88199&139185&195723&264353&10072&34&1&1&1&1&1&1&1&1&1\\
30&113&29&429&3333&11449&15630&11306&4042&700&75&10&4&1&1&1&1&1&1&1&1&1\\
31&127&32&528&4020&17304&42350&87542&158270&262703&410690&613730&885761&74677&1387&45&9&3&1&1&1&1\\
32&131&33&561&5219&21403&39349&42167&21209&5252&786&47&8&4&3&2&1&1&1&1&1&1\\
33&137&35&615&5583&23231&37766&39668&20248&5297&552&31&8&1&1&1&1&1&1&1&1&1\\
34&139&35&630&6172&25188&44702&47374&28750&9323&1961&267&32&3&2&1&1&1&1&1&1&1\\
35&149&38&697&6905&38525&138881&154706&44839&12248&1633&148&5&1&1&1&1&1&1&1&1&1\\
36&151&38&741&6743&31424&72926&81047&39200&8470&985&83&10&1&1&1&1&1&1&1&1&1\\
37&157&40&835&8805&35921&91131&125956&101503&49464&25542&3602&300&14&1&1&1&1&1&1&1&1\\
38&163&41&861&9721&65141&191575&227200&74444&21050&2439&121&6&1&1&1&1&1&1&1&1&1\\
39&167&42&903&10428&66519&198627&250421&83981&24519&3352&259&17&1&1&1&1&1&1&1&1&1\\
40&173&44&925&9765&57865&219861&618425&1493637&3138025&4160997&129892&2462&6&1&1&1&1&1&1&1&1\\
41&179&45&1035&12723&69900&173360&251966&189054&89210&52690&60940&69860&629&1&1&1&1&1&1&1&1\\
42&181&46&1091&12991&81715&369071&596594&365810&43296&6585&391&23&1&1&1&1&1&1&1&1&1\\
43&191&48&1176&15348&116722&588720&1810128&2540682&1634300&327132&9494&26&1&1&1&1&1&1&1&1&1\\
44&193&49&1279&20221&154243&642253&1811389&4020781&7619869&21743809&8337004&1142436&9690&34&4&1&1&1&1&1&1\\
45&197&50&1225&16305&118685&529537&804006&375128&97219&15641&1756&174&16&1&1&1&1&1&1&1&1\\
46&199&50&1275&15393&95200&294132&461017&353203&118576&21237&2078&152&5&1&1&1&1&1&1&1&1\\
47&211&53&1431&20471&129716&373031&639897&536759&223648&43190&4764&330&14&1&1&1&1&1&1&1&1\\
48&223&56&1596&21672&156408&590702&1604540&3479936&6581035&11350432&2243631&92958&359&2&1&1&1&1&1&1&1\\
49&227&57&1653&25293&221115&878931&1776036&2429148&4520370&7253090&770729&345&30&1&1&1&1&1&1&1&1\\
50&229&58&1739&26171&208451&1171291&2466742&2216720&491129&112817&12533&746&31&7&1&1&1&1&1&1&1\\
51&233&59&1749&26503&190383&596206&1034702&891786&389226&86298&8441&488&9&1&1&1&1&1&1&1&1\\
52&239&60&1830&29372&267698&1561616&5510739&9222901&7366619&2007701&94932&404&24&2&1&1&1&1&1&1&1\\
53&241&61&1979&38809&546491&6281521&24787263&62936401&124735849&215364097&74367005&27376541&1702999&2&1&1&1&1&1&1&1\\
54&251&63&2016&33878&264758&965348&1932458&1872501&894376&220610&24847&1301&40&1&1&1&1&1&1&1&1\\
55&257&65&2145&36089&289853&1077046&2414750&2874248&2021505&859823&151420&14679&782&23&3&1&1&1&1&1&1\\
56&263&66&2211&38824&363711&1629547&3295318&2188512&1191123&365657&55588&3891&162&7&1&1&1&1&1&1&1\\
57&269&68&2257&38405&351135&1901653&7138895&20307973&49251695&75595229&106913095&143176565&6242780&864&1&1&1&1&1&1&1\\
58&271&68&2346&38990&326872&1368328&2900734&2911574&1318748&316413&39153&2653&137&11&1&1&1&1&1&1&1\\
59&277&70&2495&44615&348973&1544371&3423489&3778726&1982886&534604&72418&5091&169&6&1&1&1&1&1&1&1\\
60&281&71&2553&46691&486703&3568771&19503671&86833859&177961341&306872855&53055757&336160&1520&6&1&1&1&1&1&1&1\\
61&283&71&2556&48048&504136&2429244&5789040&8422550&16772999&28426344&3845490&1856&51&5&1&1&1&1&1&1&1\\
62&293&74&2665&48225&468685&2688993&10354101&31703985&80145925&125119425&10957584&602241&4558&10&1&1&1&1&1&1&1\\
63&307&77&3003&60973&686435&3581683&9254016&14144918&28858373&49968934&7388646&7849&352&5&1&1&1&1&1&1&1\\
64&311&78&3081&63336&721459&5068694&24887327&94529908&297662793&813717098&1261841916&809829209&279434630&39688205&1438254&6655&39&10&3&1&1\\
65&313&79&3259&77491&872947&5134727&10511955&5808647&5078095&657099&85493&5616&292&18&5&2&1&1&1&1&1\\
66&317&80&3145&63885&689155&4302789&10447364&8733907&4501510&1218336&153083&10044&333&7&1&1&1&1&1&1&1\\
67&331&83&3486&76014&765394&3630036&9569489&12182165&7553884&2199645&300572&20811&877&44&1&1&1&1&1&1&1\\
68&337&85&3749&93685&1142877&7272661&28489765&81209509&189599989&616047989&332132160&77738924&1921708&10929&46&2&1&1&1&1&1\\
69&347&87&3828&87294&1092195&6377295&17104801&16604851&10709630&3410983&562001&44933&1711&50&1&1&1&1&1&1&1\\
70&349&88&3919&87461&995895&6928917&20600076&25074154&8378260&2572656&409481&37102&1698&88&5&1&1&1&1&1&1\\
71&353&89&3969&89921&982037&4858382&12482598&16756756&10846844&3217355&506943&38291&1277&31&1&1&1&1&1&1&1\\
72&359&90&4095&96456&1243095&9675762&44337285&99760872&110664637&46947275&4386955&37400&1470&24&1&1&1&1&1&1&1\\
73&367&92&4323&97659&1148612&7014776&27303965&77131281&178575509&362373181&119590478&10632066&131691&197&1&1&1&1&1&1&1\\
74&373&94&4471&106687&1155807&6861225&20673079&30525190&21346799&7691840&1362583&131621&6342&188&5&1&1&1&1&1&1\\
75&379&95&4560&113130&1309852&7187591&21701294&32155044&23684210&8148262&1329407&114388&5788&136&1&1&1&1&1&1&1\\
76&383&96&4656&116640&1537704&9909816&29629016&32616096&23607647&8381761&1425174&130518&6687&216&4&1&1&1&1&1&1\\
77&389&98&4717&115105&1481595&10866537&31679018&34889254&27613181&11333793&2175414&206915&9724&229&6&1&1&1&1&1&1\\
78&397&100&5095&130225&1497204&9395984&30624615&54147900&50782731&28042493&6814773&930444&67163&3128&91&2&1&1&1&1&1\\
79&401&101&5151&133197&1909099&18051659&121804835&649366397&1581410215&3078648497&5261780817&8592005051&1219157056&13838401&340&2&1&1&1&1&1\\
80&409&103&5356&159444&2984051&42350083&249145875&869391703&2195599663&4562319367&2078720651&810739411&74387757&34&3&2&1&1&1&1&1\\
81&419&105&5565&151971&1985465&12409443&41159351&66763945&53519421&21918162&4304889&412637&21137&546&19&1&1&1&1&1&1\\
82&421&106&5699&153531&2117411&17749355&95564291&293937931&664489819&1267644123&2135641251&3836551147&746908346&24243215&11339&6&1&1&1&1&1\\
83&431&108&5886&165140&2501374&22355624&116095869&299394985&385489342&197489006&24415699&300456&12995&232&5&1&1&1&1&1&1\\
84&433&109&6139&193081&2902759&22343769&65615957&58840337&55014263&15791890&2827707&259088&12816&470&17&5&1&1&1&1&1\\
85&439&110&6105&165319&2260608&15446240&53095057&89790701&71242103&28915049&5878763&634452&34968&905&17&1&1&1&1&1&1\\               
86&443&111&6216&179048&2778480&20364410&70058556&92777154&77280970&31772123&6450064&648652&28896&741&13&1&1&1&1&1&1
\end{tabular}}
\end{center}
\caption{The functions $W(p,x)$ for small values of the prime $p$. All the values in the next columns are equal to $1$.\label{eghyslist}}
\end{table}

On the assumption of the generalized Riemann hypothesis, H. L. Montgomery \cite{Montgom} proved that  $n(p)=\Omega(\log p\log\log p)$ \ie that $\overline{\lim}\ n(p)/(\log p\log\log p)>0$ and in \cite{GR}  an unconditional proof was given that
$n(p)=\Omega(\log p\log \log \log p)$. 
Combined with Lemma \ref{lowerbound} this shows that the asymptotic behavior of the function $L(p)$ is no better than $\log p\log \log \log p$.

\subsection{Experimental tests}\label{exptests}

In collaboration with Stephane Gaubert we computed the functions $W(p,x)$ for all primes up to $p=443$. We made an arborescent enumeration of the solutions
$f(2),\dots, f(L)$, for increasing values of $L$. Inadmissible sequences
were eliminated by a sieve construction. This leads to an algorithm
running in time $O(\sum_{2\leq x\leq L(p)}pxW(p,x))$.
The result is given in Table~\ref{eghyslist}.

For the ``logarithmic size"  \ie the function $\sigma(p)$,
$$
\sigma(p)=\sum_x \log W(p,x),
$$
the Gaussian approximation gives the following estimate $g(p)$, 
$$
g(p)= \frac{\left(2 \log (p)+3 \log \left(\frac{p+1}{2 p}\right)\right)^3}{12 \log ^2\left(\frac{p+1}{2 p}\right)}\simeq \frac{\left(2 \log (p)-3 \log 2\right)^3}{12 \log ^2\left(2\right)}
$$	
	
	  \begin{figure}[H]
\begin{center}
\includegraphics[scale=0.8]{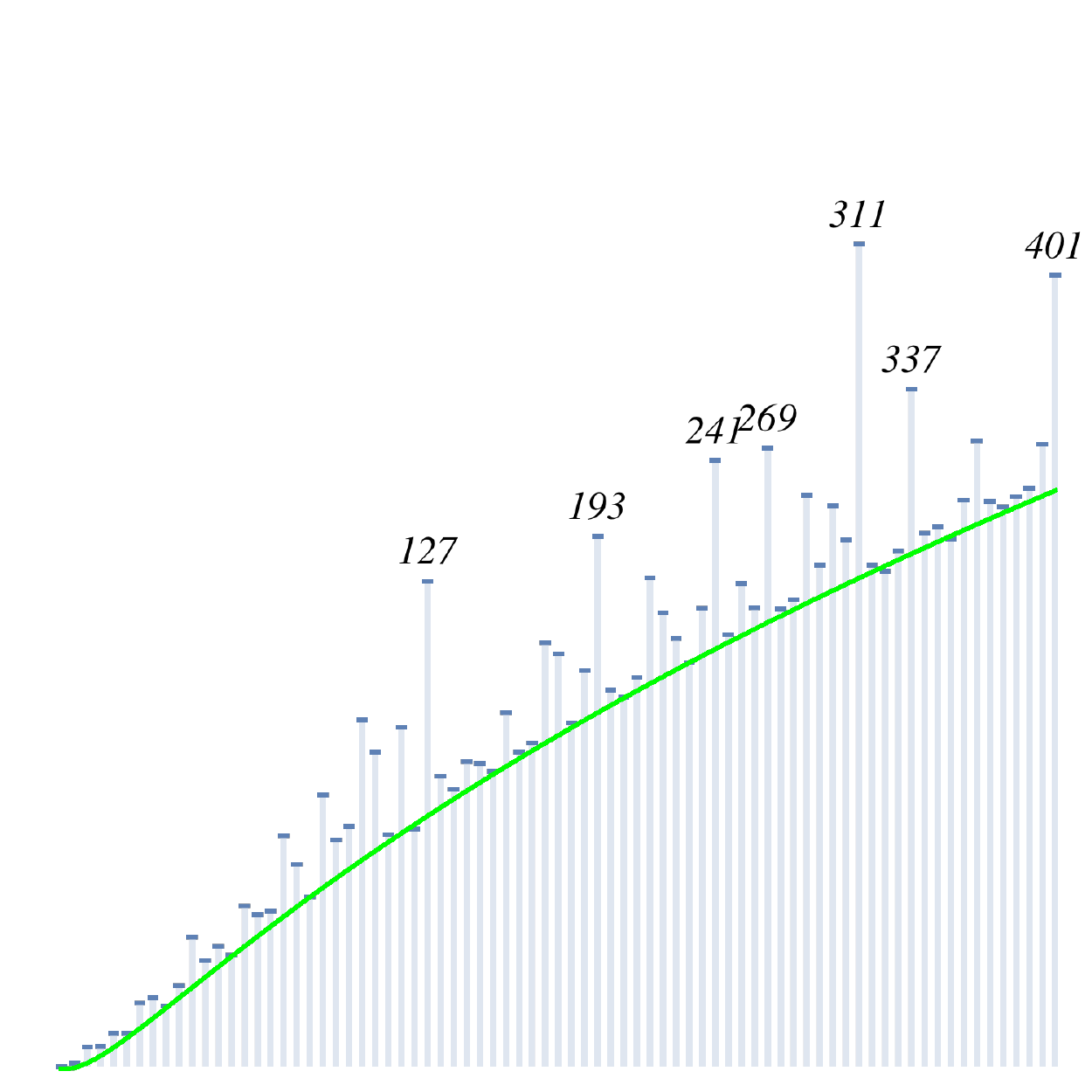}
\end{center}
\caption{Logarithmic size  for $5\leq n\leq 79$. \label{logs3} }
\end{figure}
Figure \ref{logs3} shows the graph of $g(p(n))$ in green together with the plot of $\sigma(p(n))$. It shows that there are primes such as $\{101,127,193,241,269,311,337,401\}$ for which the logarithmic size far exceeds the Gaussian estimate, but that the latter does work in many cases.  
The next Figure \ref{logsize} shows the graph of the difference $\sigma(p(n))-g(p(n))$ up to $n=79$, \ie $p(n)=401$. It is an open challenging question to understand the reason behind these non-Gaussian behaviors.
	  \begin{figure}[H]
\begin{center}
\includegraphics[scale=0.8]{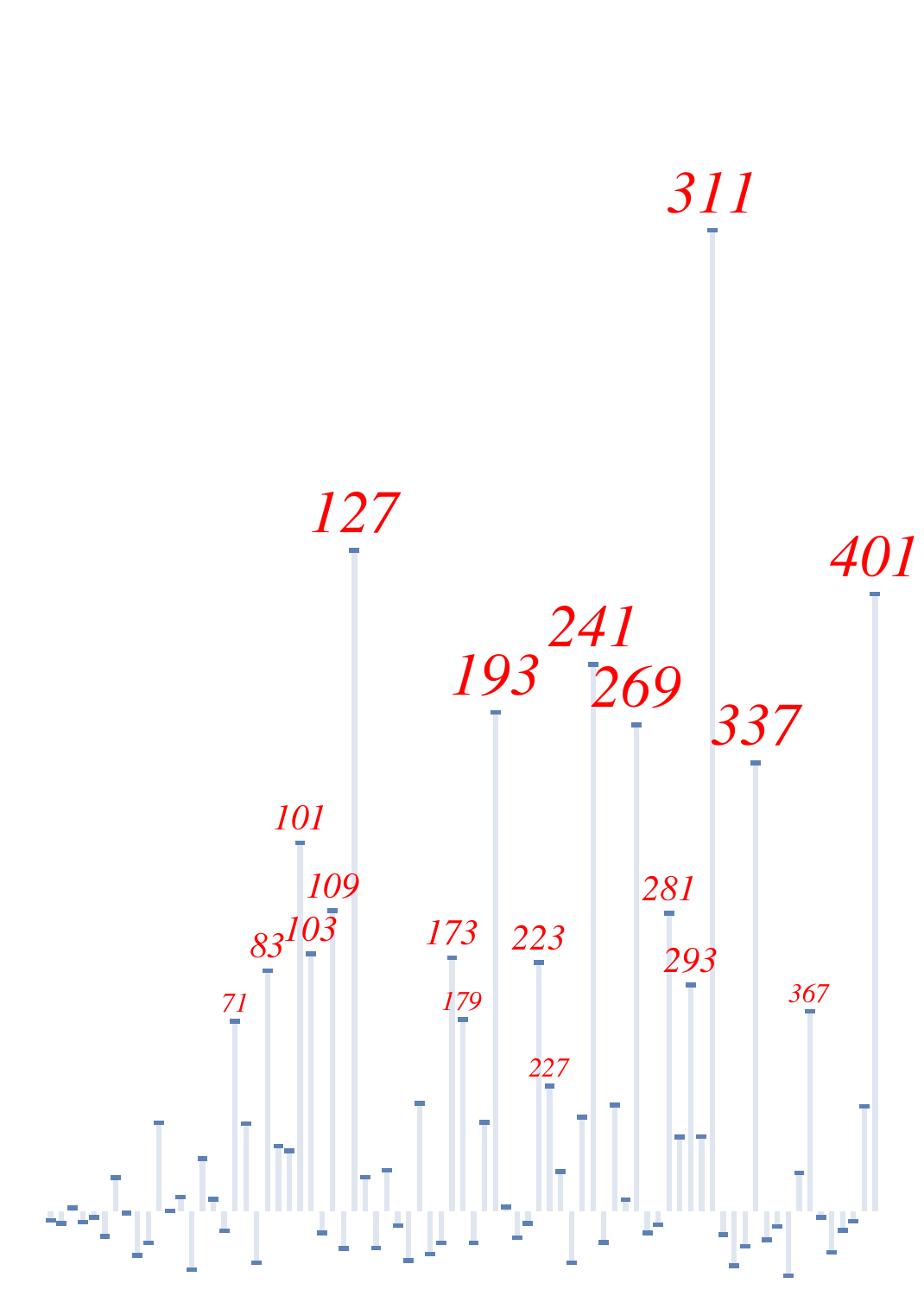}
\end{center}
\caption{Log discrepancy for $5\leq n\leq 79$. \label{logsize} }
\end{figure}

\begin{figure}[H]
\begin{center}
\includegraphics[scale=0.9]{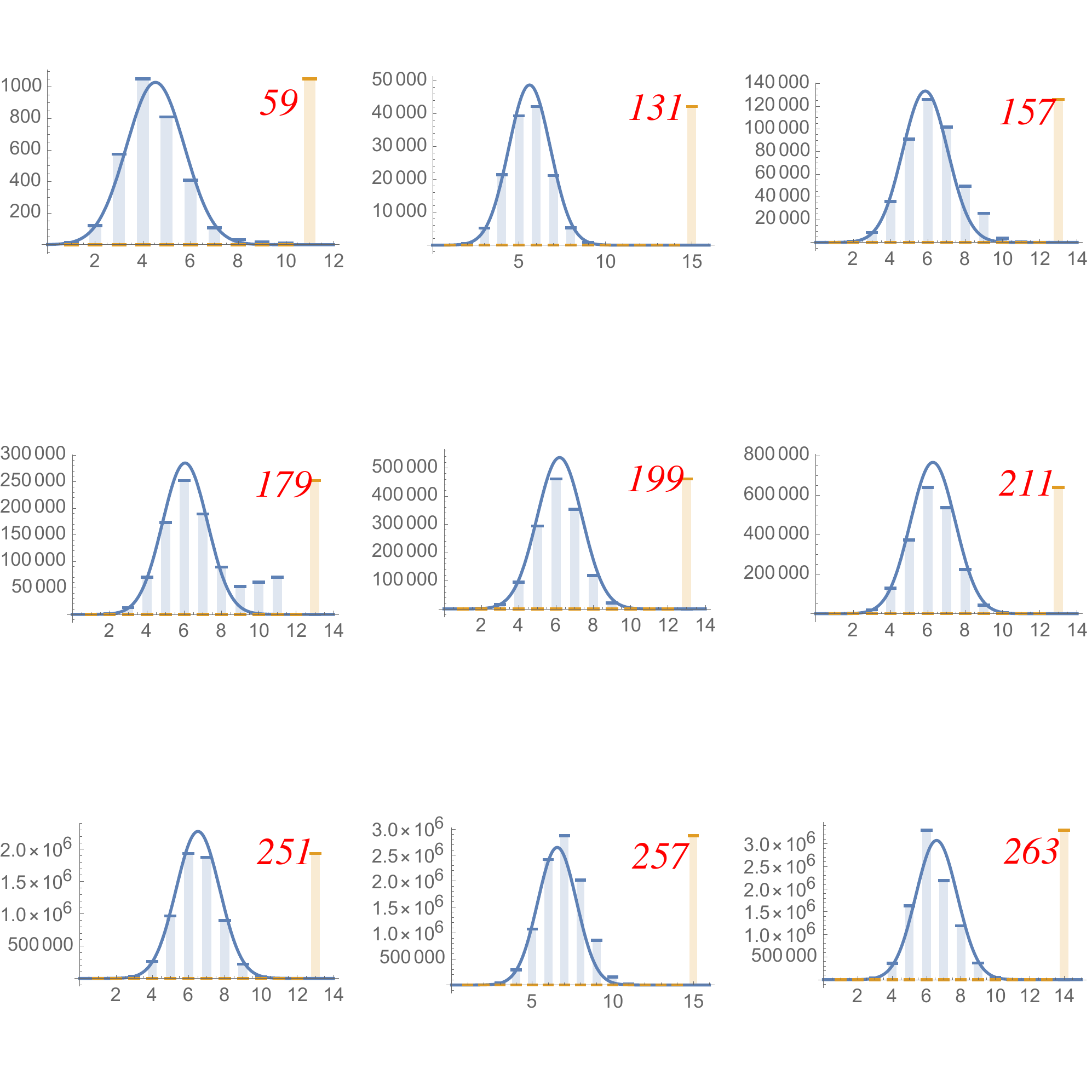}
\end{center}
\caption{Primes with good gaussian approximation. Graphs of the functions $W(p,x)$ and $G(p,x)$ with $x$ a continuous variable for the Gaussian. \label{ronkin2} }
\end{figure}

\begin{figure}[H]
\begin{center}
\includegraphics[scale=0.7]{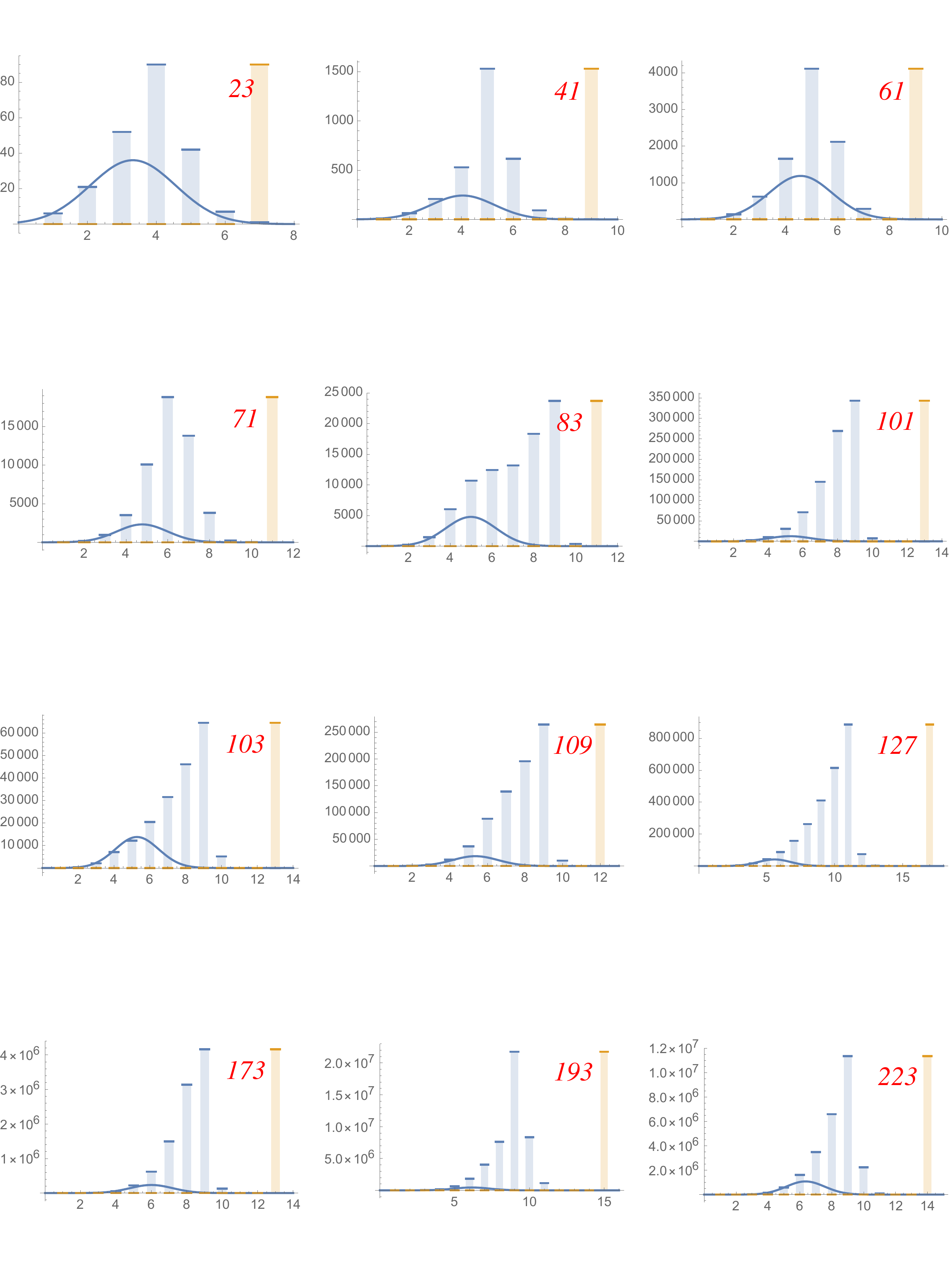}
\end{center}
\caption{Primes with bad gaussian approximation. \label{ronkin2} }
\end{figure}


\begin{thebibliography}{99}

\bibitem{Carlitz} L.~Carlitz, {\em A theorem on permutations in a finite field}, Proc. Amer. Math. Soc. 11 (1960), 456--459. Errata ibid. 999--1000.
	 	
	 	\bibitem{Jones} G.Jones, {\em Paley and the Paley graphs}, ArXiv 1702.00285
 
 \bibitem{GR} S.~Graham, C.~Ringrose {\em Lower bounds for least quadratic nonresidues}. 
 Analytic number theory (Allerton Park, IL, 1989), 269--309, 
Progr. Math., 85, Birkhäuser Boston, Boston, MA, 1990. 


\bibitem{Montgom} H. L. Montgomery, {\em Topics in multiplicative number theory}, Lecture Notes in Math., 227, Springer, Berlin, 1971.
\end{thebibliography}
\end{document}